\input amstex
\documentstyle{amsppt}
\magnification=\magstep1
 \hsize 13cm \vsize 18.35cm \pageno=1
\loadbold \loadmsam
    \loadmsbm
    \UseAMSsymbols
\topmatter
\NoRunningHeads
\title A note on $p$-adic $q$-integral associated with $q$-Euler numbers
\endtitle
\author
  Taekyun Kim
\endauthor

\abstract We show that $q$-Euler numbers can be represented as an
integral by the $q$-analogue of the ordinary $p$-adic fermionic
measure, whence we give an answer to the question of readers  to
ask us for the fermionic $p$-adic $q$-integral on $\Bbb Z_p$.
Sometimes several readers give the nonsense comments to us. But we
do not write their names in this paper. I would like to tell the
readers that papers are different from books. That is, the readers
are need their efforts to understand the contents of paper.

\endabstract
\endtopmatter

\document

{\bf\centerline {\S 1. Introduction}}

 \vskip 20pt

Throughout  this paper $\Bbb Z$, $\Bbb Q$, $\Bbb Z_p$, $\Bbb Q_p$
and $\Bbb C_p$ will respectively denote the ring of rational
integers, the field of rational numbers, the ring $p$-adic
rational integers, the field of $p$-adic rational numbers and the
completion of the algebraic closure of $\Bbb Q_p$. Let $v_p$ be
the normalized exponential valuation of $\Bbb C_p$ such that
$|p|_p=p^{-v_p(p)}=p^{-1}$. If $q\in\Bbb C_p$, we normally assume
$|q-1|_p<p^{-\frac{1}{p-1}},$ so that $q^x =\exp (x\log q)$ for
$|x|_p\leq 1$. We use the notation
$$[x]_q=\frac{1-q^x}{1-q}, \text{ and  }
[x]_{-q}=\frac{1-(-q)^x}{1+q}.$$ Hence, $\lim_{q\rightarrow
1}[x]_q =1, $ for any $x$ with $|x|_p\leq 1$ in the present
$p$-adic case.

Let $p$ be a fixed odd prime. For $d(= odd)$ a fixed positive
integer with $(p,d)=1$, let
$$\split
& X=X_d = \lim_{\overleftarrow{N} } \Bbb Z/ dp^N \Bbb Z ,\cr & \
X_1 = \Bbb Z_p , \cr  & X^\ast = \underset {{0<a<d p}\atop
{(a,p)=1}}\to {\cup} (a+ dp \Bbb Z_p ), \cr & a+d p^N \Bbb Z_p =\{
x\in X | x \equiv a \pmod{dp^n}\},\endsplit$$ where $a\in \Bbb Z$
lies in $0\leq a < d p^N$.

In this paper we prove  that $$\mu_{-q}(a+dp^N \Bbb
Z_p)=(1+q)\frac{(-1)^aq^a}{1+q^{dp^N}}
=\frac{(-q)^a}{[dp^N]_{-q}},$$ is distribution on $X$ for $q\in
\Bbb C_p$ with $|1-q|_p< p^{-\frac{1}{p-1}}. $ This distribution
yields an integral as follows:
$$I_{-q}(f)=\int_{\Bbb Z_p} f(x)d\mu_{-q}(x)=\lim_{N\rightarrow
\infty}\frac{1}{[p^N]_{-q}}\sum_{x=0}^{p^N-1}f(x)(-q)^x, \text{
for $f\in UD(\Bbb Z_p)$ },
$$
which has a sense as we see readily that the limit is convergent.

For $q=1$,  we have fermionic $p$-adic integral on $\Bbb Z_p$ as
follows:
$$I_{-1}=\int_{\Bbb
Z_p}f(x)d\mu_{-1}(x)=\lim_{N\rightarrow
\infty}\sum_{x=0}^{p^N-1}f(x)(-1)^x. $$ In view of notation,
$I_{-1}$ can be written symbolically as
$I_{-1}(f)=\lim_{q\rightarrow -1} I_{q}(f).$ Finally, we will
introduce $q$-extension of Euler numbers by using $I_{-q}(f)$

\vskip 20pt

{\bf\centerline {\S 2. A note on fermionic  $p$-adic $q$-integral
on $\Bbb Z_p$ }}

 \vskip 20pt

For any positive integer $N$, $q\in\Bbb C_p$ with
$|1-q|_p<p^{-\frac{1}{p-1}}, $ we set
$$\mu_{-q}(a +dp^N\Bbb
Z_p)=(1+q)\frac{(-1)^aq^a}{1+q^{dp^N}}=\frac{(-q)^a}{[dp^N]_{-q}},
\tag *$$ and this can be extended  to a distribution on $X$.

We show that $\mu_{-q}$ is distribution on $X$.  For this it
suffices to check that
$$\sum_{i=0}^{p-1}\mu_{-q}(a+idp^N+dp^{N+1}\Bbb
Z_p)=\mu_{-q}(a+dp^N\Bbb Z_p). $$ The left hand side is equal to
$$(1+q)\sum_{i=0}^{p-1}\frac{(-1)^{a+idp^N}}{1+q^{dp^{N+1}}}q^{a+idp^N}=(1+q)\frac{(-1)^aq^a}{1+q^{dp^{N+1}}}
\sum_{i=0}^{p-1}q^{idp^N}(-1)^i. \tag 1$$ Since we see that
$$\frac{1}{1+q^{dp^{N+1}}}=\frac{1}{1+q^{dp^N}}\cdot\frac{1+q^{dp^N}}{1+q^{dp^{N+1}}}.\tag
2$$

From (1) and (2) we derive
$$\eqalignno{
&\sum_{i=0}^{p-1}\mu_{-q}(a+idp^N+dp^{N+1}\Bbb
Z_p)=(1+q)\frac{(-1)^aq^a}{1+q^{dp^{N+1}}}\sum_{i=0}^{p-1}q^{idp^N}(-1)^i\cr
&=(1+q)\frac{(-1)^aq^a}{1+q^{dp^N}}\cdot\frac{1+q^{dp^N}}{1+q^{dp^{N+1}}}\sum_{i=0}^{p-1}(-1)^iq^{idp^N}
=\frac{(-1)^aq^a(1+q)}{1+q^{dp^N}}=\mu_{-q}(a+dp^N\Bbb Z_p ).}
$$
This distribution yields an integral for each non-negative integer
$m$ in the case $d=1$,
$$\int_{\Bbb Z_p}[a]_q^md\mu_{-q}(a) =\lim_{N\rightarrow \infty}\sum_{a=0}^{p^N-1}[a]_q^m\frac{(-q)^a(1+q)}{1+q^{p^N}}
=I_{-q}([a]_q^m), \tag3$$ which has a sense as we see readily that
the limit is convergent.

Also, we easily see that (*) is  distribution on $X$ for $q\in\Bbb
C_p$ with $|1-q|_p <1.$ We now define a $q$-Euler numbers
$E_{m,q}\in\Bbb C_p$ by making use of this integral:
$$ I_{-q}([a]_q^m)=\int_{\Bbb Z_p}[a]_q^md\mu_{-q}(a)=E_{m,q}. $$
Note that $\lim_{q\rightarrow 1}E_{m,q}=E_m,$ where $E_m$ are the
$m$-th ordinary Euler numbers.

The generating function $F_q(t)$ of $E_{k,q}$,
$$F_q(t)=\sum_{k=0}^{\infty}E_{k,q}\frac{t^k}{k!}, \tag4$$
is given by
$$F_q(t)=\lim_{\rho\rightarrow\infty}\frac{1}{[p^{\rho}]_{-q}}\sum_{i=0}^{p^{\rho}-1}(-q)^i e^{[i]_qt}, \tag5$$
which satisfies the $q$-difference equation
$$F_q(t)=-qe^tF_q(qt)+1. \tag6$$
If $q=1$ in Eq.(3)-(6), then we have
$$F_q(t)=\frac{2}{e^t+1}=\sum_{n=0}^{\infty}E_n\frac{t^n}{n!}. $$
Let $\chi$ be a primitive Dirichlet character with conductor
$d(=odd)\in\Bbb Z_{+}$, the set of natural numbers. Then we also
define  a generalized $q$-Euler numbers $E_{n,\chi, q}$ as
$$E_{m,\chi,
q}=\int_{X}\chi(a)[a]_q^md\mu_{-q}(a)=\lim_{N\rightarrow\infty}\sum_{a=0}^{dp^N-1}[a]_q^m\chi(a)\frac{(-q)^a}
{[dp^N]_{-q}},$$ where $E_{m,\chi, q}$ are the $m$-th generalized
$q$-Euler numbers attached to $\chi$.

The $q$-Euler polynomials in the variable $x$ in $\Bbb C_p$ with
$|x|_p\leq 1$ are defined by
$$\int_{\Bbb Z_p} [x+t]_q^n d \mu_{-q} (t) = E_{n,q}(x) . $$
These can be written as
$$E_{n,q}(x)=\sum_{l=0}^n \binom nl   q^{lx}E_{l,q}[x]_q^{n-l}. $$
Indeed we see
$$\int_{\Bbb Z_p}[x+t]_q^nd\mu_{-q}(t)=\sum_{k=0}^n\binom nk
[x]_q^{n-k}q^{kx}\int_{\Bbb Z_p}[t]_q^kd\mu_{-q}(t).$$ For the
integral $I_{-q}$ we first see:

\proclaim{Theorem 1} For $m \geq 0$, $q\in \Bbb C_p$ with
$|1-q|_p<p^{-\frac{1}{p-1}}, $ we have

$$E_{m,q}=\int_{\Bbb Z_p}[x]_q^md\mu_{-q}(x)=[2]_q
\left(\frac{1}{1-q}\right)^m\sum_{k=0}^m\binom mk
(-1)^k\frac{1}{1+q^{k+1}}.$$

Proof. We see
$$\split
&\frac{1+q}{1+q^{p^N}}\sum_{a=0}^{p^N-1}[a]_q^mq^a(-1)^a=\frac{1+q}{1+q^{p^N}}\frac{1}{(1-q)^m}
\sum_{a=0}^{p^N-1}(-1)^aq^a\sum_{j=0}^m \binom mj (-1)^jq^{aj}\cr
  &=\frac{1}{(1-q)^m}\frac{1+q}{1+q^{p^N}}\sum_{j=0}^m\binom
  mj(-1)^j\frac{1+q^{(j+1)p^N}}{1+q^{j+1}}.
\endsplit$$
Since $\lim_{N\rightarrow \infty}q^{p^N}=1$ for
$|1-q|_p<p^{-\frac{1}{p-1}}, $ our assertion follows.
\endproclaim

\proclaim{Lemma 2} For $q\in\Bbb C_p$ with $|1-q|_p<1$, we have
$\lim_{N\rightarrow \infty}q^{p^N}=1.$

Since $$ q^{p^N}=(q-1+1)^{p^N}=\sum_{l=0}^{p^N}\binom
{p^N}{l}(q-1)^l. $$
\endproclaim

From Theorem 1 and Lemma 2, we derive the following:

\proclaim{ Corollary 3} For $q\in\Bbb C_p$ with $|1-q|_p<1$, $m\in
\Bbb Z_{+},$ we have
$$\int_{\Bbb Z_p}[x]_q^md\mu_{-q}(x)=[2]_q
\left(\frac{1}{1-q}\right)^m\sum_{k=0}^m\binom mk
(-1)^k\frac{1}{1+q^{k+1}}.$$

\endproclaim

For $f \in UD(\Bbb Z_p)$, let us start with expression
$$\frac{1}{[p^N]_{-q}}\sum_{0\leq j<p^N}(-q)^jf(j)=\sum_{0\leq
j<p^N}f(j)\mu_{-q}(j+p^N\Bbb Z_p) ,$$ representing $q$-analogue of
Riemann sums for $f$.

The fermionic $p$-adic $q$-integral on $\Bbb Z_p$ will be defined
as limit ( $N\rightarrow \infty$) of these sums, when it is
exists. A fermionic $p$-adic $q$-integral of function $f \in
UD(\Bbb Z_p)$ on $\Bbb Z_p$ is defined as
$$ I_{-q}(f)=\int_{\Bbb Z_p}f(x)d\mu_{-q}(x)=\lim_{N\rightarrow
\infty}\frac{1}{[p^N]_{-q}}\sum_{x=0}^{p^N-1}f(x0q^x(-1)^x.\tag 7
$$
Note that if $f_n \rightarrow f$ in $UD(\Bbb Z_p)$; then
$$\int_{\Bbb Z_p} f_n d\mu_{-q}(x)\rightarrow \int_{\Bbb
Z_p}f(x)d\mu_{-q}(x). $$

From (7), we derive the following theorem:
 \proclaim{ Theorem 2}
For $f \in UD(\Bbb Z_p)$, we have
$$qI_{-q}(f_1 )+I_{-q}(f)=[2]_qf(0), $$
where $f_1(x)$ is translation with $f_1(x)=f(x+1).$
\endproclaim

Remark. In [1], $p$-adic $q$-integral on $\Bbb Z_p$ is defined by
$$I_q(f)=\int_{\Bbb Z_p}f(x)d\mu_{q}(x)=\lim_{N\rightarrow
\infty}\frac{1}{[p^N]_q}\sum_{x=0}^{\infty}f(x)q^x .$$ For $q=1$
in Eq. (7), we have
$$I_{-1}(f)=\int_{\Bbb Z_p}f(x)d\mu_{-1}(x)=\lim_{N\rightarrow
\infty}\sum_{x=0}^{p^N-1}f(x)(-1)^x. $$ In view of notation,
$I_{-1}(f)$ can be written symbolically as
$I_{-1}(f)=\lim_{q\rightarrow -1} I_q(f). $

\vskip 20pt

\quad Taekyun Kim

\quad EECS, Kyungpook National University, Taegu 702-701, S. Korea

\quad e-mail:\text{ tkim$\@$knu.ac.kr; tkim64$\@$hanmail.net}

\enddocument